\documentclass[12pt,a4paper]{article}
\usepackage{titlesec}
\titleformat{\section}[block]{\Large\bfseries\filcenter}{\thesection}{1em}{}
\titleformat{\subsection}[hang]{\bfseries\filright}{\thesubsection}{1em}{}

\usepackage{amsfonts,amssymb,amsmath,amsthm,textcomp,mathrsfs}

\usepackage{graphicx}
\usepackage{hyperref}

\titlelabel{\thetitle.\quad}

\hypersetup{
    colorlinks,
    linkcolor={blue!50!black},
    citecolor={blue!50!black},
    urlcolor={blue!80!black}
}
\usepackage{longtable}
\usepackage{csquotes}
\usepackage{color,xcolor}
\definecolor{linkclr}{rgb}{0,0,0.8}
\definecolor{citeclr}{rgb}{0,0.5,0}
\definecolor{urlclr}{rgb}{0.8,0.4,0.4}
\definecolor{Plum}{HTML}{8E4585}
\definecolor{Mahogany}{HTML}{C04000}
\definecolor{Maroon}{HTML}{800000}
\definecolor{BrickRed}{HTML}{841F27}

\renewcommand{\thesection}{\textcolor{black}{\arabic{section}}}
\renewcommand{\thesubsection}{\textcolor{black}{\arabic{section}.\arabic{subsection}}}

\usepackage[PetersLenny]{fncychap} 
\ChNameVar{\fontsize{18}{12}\usefont{OT1}{phv}{m}{n}\selectfont}
\ChNumVar{\fontsize{60}{62}\usefont{OT1}{ptm}{m}{n}\selectfont\hspace{2pt}}
\ChTitleVar{\centering\Huge\sf}\ChRuleWidth{0.5pt}\ChNameUpperCase

\usepackage{geometry}
\newtheorem{thm}{Theorem}[section]

\newtheorem{lem}[thm]{Lemma}
\newtheorem{prop}[thm]{Proposition}

\theoremstyle{definition}

\addtocounter{section}{0}

 2

\usepackage{setspace}
\usepackage{appendix}

\newcommand{\Address}{{
  
\bigskip
\footnotesize
\noindent
\textsc{Harish-Chandra Research Institute, HBNI, Chhatnag Road, Jhunsi, Allahabad -211 019, India.}\par\nopagebreak
  \textit{email }: \texttt{mallesham@hri.res.in}
  }}

 \newcommand\blfootnote[1]{%
   \begingroup
   \renewcommand\thefootnote{}\footnote{#1}
   \endgroup
 }

\begin{document}
\singlespacing

\begin{center}
{\Large {\bf Primes in Sumsets}} \\
\vspace{1mm}
{\em Kummari Mallesham}
\end{center}

\vspace{5mm} 
\begin{abstract}
\noindent 
We obtain an upper bound for the number of pairs $ (a,b) \in {A\times B} $ such that $ a+b $ is a prime number, where $ A, B \subseteq \{1,...,N \}$ with $|A||B| \, \gg \frac{N^2}{(\log {N})^2}$, $\, N \geq 1$ an integer. This improves on a bound given by Balog, Rivat and S\'ark\"ozy.

\end{abstract}

\blfootnote{2010 {\it Mathematics subject classification}: Primary 11N36; Secondary 11N13 }
\blfootnote{ {\it Key words and Phrases}  : Large sieve inequality, Primes, Brun-Titchmarsh inequality }

\section{Introduction}
\label{intro}

\noindent
Let  $ A, B$ be subsets of $\{1,...,N \}, \, N \geq 1$ an integer, and let us write  $ P_{N;A,B} $  for the number of pairs $ (a,b) $ in $ A \times B $ such that $ a+b $ is a prime number. One might expect that $ P_{N;A,B} $ is about $ \frac{|A| |B|}{\log{N}} $, but there are subsets $ A, B $ of $ \{1,...,N \} $ such that $ P_{N;A,B} $ is larger than what one expects. In \cite[Section~5]{BRS}, A. Balog, J. Rivat and A. S\'ark\"ozy treated the problem of bounding $ P_{N;A,B} $ from above. Indeed, let us set $R$ to be the quantity $\frac{1000N}{|A|^{1/2}|B|^{1/2}}$. Then the following theorem is proved as Theorem 6 in \cite{BRS} using the linear sieve.

\begin{thm}[Balog, Rivat and S\'ark\"ozy]
\label{bsrthm}
We have  $ P_{N;A,B} \ll  \frac{ |A||B| } {\log N} \, R $.
\end{thm}

\noindent
We begin by observing that this bound can be obtained using a very elementary counting argument, which we give in Subsection \ref{brsproof} below. Also note that this bound is implied by the trivial estimate $P_{N;A,B} \leq |A||B|$ unless $|A||B| \gg N^2 / (\log N)^2$. In Section \ref{mainproof} we apply a method of D.S. Ramana and O. Ramar\'e \cite{SR}, which was originally used to obtain an upper  bound for additive energy of dense subsets of primes, to prove following theorem, which is our main result.

\begin{thm} 
\label{mall1}
Let $A,B$ be subsets of $\{1,\ldots,N\}$, where $N \geq 1$ is an integer and suppose further that $|A||B| \, \gg \frac{N^2}{(\log {N})^2}$. Then we have that 
$ P_{N;A,B} \ll \frac{ |A||B| } {\log N} \log \log {R}.$ 
\end{thm}

\vspace{2mm}
\noindent
This upper bound for $P_{N;A,B}$ is optimal, up to the implied constant, in general. In fact, we have the following proposition, which corrects the conclusion of Example 2, page 36 of \cite{BRS}.

\begin{prop}
Let $N$ be a positive integer, $k \ll \log \log{N}$ be an integer and $m_{k} = \prod_{p \leq k} p$. Then if $A=\{1 \leq a \leq N : a \equiv 0 \, {\rm mod} \, m_{k}\}$ and $ B=\{1 \leq b \leq N : b \equiv 1 \, {\rm mod} \, m_{k}\}$, we have that $$ P_{N;A,B} \gg \frac{ |A||B| } {\log N} \log \log{R} . $$
\end{prop}

\vspace{2mm}
\noindent
{\sc Proof.---} If $r_{A,B}(n)$ is the number of pairs $(a,b) \in A \times B$ such that $a+b=n$ then

\begin{equation}
\label{o}
P_{N;A,B} \geq \sum_{\substack{\frac{N}{2} \leq p \leq N \\ p \, \equiv \, 1 \, {\rm mod} \, m_{k}}} r_{A,B}(p).
\end{equation}

\noindent
We observe that for any integer $n \equiv 1\, \text{mod} \,m_k$ and $n \geq 2m_k$,  we have $$ r_{A,B}(n) \geq [\frac{n}{m_{k}}] \geq \frac{n}{2 m_{k}}.$$ Using this lower bound for $r_{A,B}(p)$ in \eqref{o} when $N \geq 4m_k$  we get 

\begin{equation}
\label{tp}
P_{N;A,B} \geq \frac{1}{2 m_{k}} \,\sum_{\substack{\frac{N}{2} \leq p \leq N, \\ p \, \equiv \, 1 \, {\rm mod} \, m_{k}.}} p \geq  \frac{N}{4 m_{k}} \,\sum_{\substack{\frac{N}{2} \leq p \leq N, \\ p \, \equiv \, 1 \, {\rm mod} \, m_{k}.}} 1.
\end{equation}

\noindent
By the Chebyshev bound $\log m_k = \sum_{p \leq k} \log p \ll \log\log N$. Thus on using the Siegel-Walfisz theorem (see \cite[page 419]{IK}),  we have 
 
\begin{equation} 
\label{t}
\sum_{\substack{\frac{N}{2} \leq p \leq N, \\ p \, \equiv \, 1 \, {\rm mod} \, m_{k}.}} 1 \gg \frac{N}{ \phi(m_{k}) \, \log{N}}.
\end{equation}

\noindent
Merten's formula gives the the upper bound $\phi(m_{k}) \ll \frac{m_{k}}{\log \log{m_{k}}} $. Also, $|A| \, \sim \frac{N}{m_{k}}, |B| \, \sim \frac{N}{m_{k}}$ and therefore $R \sim m_k$, from the definition of $R$. From  \eqref{tp} and \eqref{t} we then get  

$$P_{N;A,B} \gg \frac{N^2}{\phi(m_k) m_k \log N} \gg \frac{ |A||B| } {\log N} \log \log {R}, $$ which proves the proposition.

\vspace{2mm}
\noindent
In Section 2 we record some preliminaries, mainly taken from \cite{SR}, and give our proof of Theorem \ref{bsrthm}.  In Section 3 we reduce the proof of Theorem \ref{mall1} to the case of subsets which are well distributed to certain moduli. We then use this reduction to complete the proof of the Theorem \ref{mall1} in Section 4, which is our final section.

\vspace{2mm}
\noindent
Throughout this article we use $e(z)$ to denote $e^{2\pi i z}$ for any complex number~ $z$. Further, all constants implied by the symbols $\ll$, $\gg$ and the $O$ notation are absolute except when dependencies are indicated, either in words or by subscripts to these symbols. The Fourier transform of the characteristic function of a subset $A$ of $\mathbb{Z}$ is denoted by $\widehat{A}$ and is defined by $\widehat{A}(t) = \sum_{a \in A} e(at)$.  Finally, the notations $[a,b], (a,b]$ etc. will denote intervals in ${\mathbb Z}$  with end points $a$, $b$ unless otherwise specified.

\section{Preliminaries}

\label{preli}

\subsection{The large sieve inequality}\label{large}

The following is the classical large sieve inequality, which is proved on \cite[page 68]{Tennen}, for example.

\vspace{2mm}
\noindent
Let $N \geq 1 $ be a integer and $Q \geq 1 $ be a real number. Then for any sequence of complex numbers  
$\lbrace a_{n}\rbrace _{n=1}^{N}$ and real number $\alpha$ if we set  $$ S(\alpha) = \sum_{\substack{1 \leq n \leq N}} a_{n} \, e(n \alpha),$$ we have 
 
\begin{equation}
\label{largeineq}
\sum_{\substack{1 \leq q \leq Q}}  \, \, \sum_{\substack{a \, {\rm mod^{*}} \, q}} |S(a/q)|^{2}  \, \leq (N + Q^2) \sum_{\substack{1 \leq n \leq N}} |a_{n}|^2  .
\end{equation}

\subsection{The Brun-Titchmarsh inequality} \label{siegel}

\noindent
If $ q,a $ are positive integers with $(a,q)=1 $, then for all $ q \leq x $, we have
  
\begin{equation}
\pi(x;q,a) \leq \frac{2x}{\phi(q) \, \log{(x/q)}} \,, 
\end{equation}

\noindent
where $\pi(x;q,a)$ denotes the number of primes not exceeding $x$ and congruent to $a$ modulo $q$. For a proof see \cite[page 121]{Mont}. In particular, we have 

\begin{equation}
\label{siegel} 
\pi(x;q,a) \leq \frac{4x}{\phi(q) \, \log {x} }  \,,
\end{equation}
when $ q \leq x^{\frac{1}{2}}$.

\subsection{An arithmetical function } \label{asymp}

For any integer $q \geq 1$ and a positive real number $L \geq 1$,  let us set 

\begin{equation}
\label{arith}
\omega (q, L) = - \sum_{\substack{1 \leq l \leq L, \\ l \, \equiv \, 0 \, {\rm mod } \, q }} \frac{\mu(l) \log{l}}{l} \; ,
\end{equation}

\noindent
where $\mu$ is the M{\"o}bius function. We then have the following estimates for $ \omega (q,L)$, proved in \cite[Section~2.1] {SR}. Here $\nu(q)$ denotes the number of prime divisors of $q$.

\begin{lem}
\label{asymptotic}

(i)~ For $ 1 \leq q \leq L^{1/2}$ , we have the asymptotic formula

\begin{equation}
\label{asymptotic1}
\omega(q,L)= \frac{\mu(q)}{ \phi (q)} + O_{\alpha} \left( \frac{2^{\nu(q)} \log{2q}}{q \, (\log{L})^{\alpha}}\right) \, ,
\end{equation}

\noindent
for any $\alpha >0$ and (ii)~ for any $q , L \geq 1 $ , we have

\begin{equation}
\label{asymptotic2}
| \,\omega (q, L) \, |  \, \leq \, \frac{(\log{2L})^{2}}{q} \, .
\end{equation} 

\end{lem}

\subsection{ An application of Davenport's bound} \label{upperbo}
Let $N \geq 1 $ be an integer, $L=N^{1/2}$ and set $\Lambda^{\flat}(n) = - \sum_{\substack{ d|n,\\ d > L.}}  \mu  (d) \log{d} $ for any integer $n$.

\noindent
The following lemma gives a uniform bound for the Fourier transform of the restriction of $n \mapsto \Lambda^{\flat}(n)$ to the interval $[1,2N]$.

\begin{lem}
\label{ftbound}
We have $$ \sum_{\substack{1 \leq n \leq 2N}} \Lambda^{\flat}(n) \, e(nt)  \, \ll  \, \frac{N}{(\log N)^{100}} , $$ for all $t \in [0,1]$.
\end{lem} 

\vspace{2mm}
\noindent
{\sc Proof.---} The lemma follows from Davenport's classical bound for $\sum_{1 \leq n \leq x} \mu(n)e(nt)$, given by Theorem 13.10 on page 348 of \cite{IK}, by an integration by parts. See \cite[Section~3]{SR} for the details.
  
\subsection{ An optimisation principle} \label{opt}

In our proof of the Theorem \ref{mall1} we will use an optimization principle, a minor variant on a similar principle from \cite{SR}. We state this principle with the aid of the following notation.

\vspace{2mm}
\noindent
Suppose that $n, m \geq 1$ are integers and let $P_{1}, P_{2}, D_{1} $ and  $ D_{2}$  be real numbers $> 0$. Further let 
\begin{equation*}
\mathcal{K}_{1}= \left \lbrace (x_{1}, \ldots,x_{m}) \in \mathbb{R}^{m}  \,  : \, \sum_{i=1}^{m} x_{i} = P_{1}, \, \, 0 \leq x_{i} \leq D_{1} \, \, \text{for all $i$} \right \rbrace \, ,
\end{equation*}

\noindent
and
\begin{equation*}
\mathcal{K}_{2}= \left \lbrace (x_{1}, \ldots,x_{n}) \in \mathbb{R}^{n}  \,  : \, \sum_{i=1}^{n} x_{i} = P_{2}, \, \, 0 \leq x_{i} \leq D_{2} \, \, \text{for all $i$} \right \rbrace \, .
\end{equation*}

\vspace{2mm}
\noindent
Let us also assume that $\mathcal{K}_{1}$ and $\mathcal{K}_{2}$ are non-empty sets. Then $\mathcal{K}_{1}, \mathcal{K}_{2}$ are compact and convex subsets of $\mathbb{R}^{m}, \, \mathbb{R}^{n}$ respectively. Then we have :

\begin{lem} \label{bil}
If $f: \mathbb{R}^m \times \mathbb{R}^n \mapsto \mathbb{R}$ a bilinear form with real coefficients $\alpha_{ij}$ defined by  $f(x,y) =\sum_{\substack{1 \leq i \leq m , \\ 1 \leq j \leq n}} \alpha_{ij} x_i y_j$ then 

\vspace{1mm}
\noindent
$(i)$ there are extreme 
points $x^{*}$ and $y^{*}$ of ${\mathcal{K}_{1}}$ and $\mathcal{K}_{2}$ respectively  so that $f(x,y) \leq 
f(x^{*},y^{*})$ for all $x \in {\mathcal {K}_{1}} , \, y \in {\mathcal{K}_{2}}$.

\vspace{1mm}
\noindent
$(ii)$ If $x^{*} = (x_1^{*}, x_2^{*}, \ldots, x_m^{*})$ is an extreme point 
of ${\mathcal {K}_{1}}$ then, excepting at most one ~$i$, we have either $x_i^{*} = 0$ or $x_i^{*} = D_{1}$ for each $i$. Also, if $ l $ is the 
number of $i$ such $x_i^{*} \neq 0$ then $l D_{1} \geq P_{1} > (l-1) D_{1}$. A similar result holds for the extreme points of $\mathcal{K}_{2}$.
 
\end{lem}

\vspace{2mm}
\noindent
{\sc Proof.---} An easy modification of the proof of Proposition 2.2 of \cite{SR}.

\subsection{A local problem} \label{localprob}
  
\noindent
Let $R \geq 1000$ be real number, which we will eventually take to be $\frac{1000N}{|A|^{1/2}|B|^{1/2}}$, as in Section 1. We then set  
  
  \begin{equation}
  \label{umq}
  U= \prod_{p \leq R} p,\;\; \, M_{R}= \left(\frac{R \, \log R}{\log \log R}\right)^2 \; \; \text{and}\;\;  Q = \log R \, \log \log R. 
  \end{equation}
 
 \noindent 
 Also, let $I$ be the set of prime numbers not exceeding $R$ and  let $J$ be a given subset of $I$. We then write $T_{J}(\mathcal{X}, \mathcal{Y})$ for the number of pairs $(x,y)$ in $\mathcal{X} \times \mathcal{Y}$ such that $x+y \not\equiv 0 \, (\text{mod} \, p)$ for each $p$ in $J$. The following lemma, which gives an upper bound for $T_{J}(\mathcal{X}, \mathcal{Y})$, is Proposition 2.3 of \cite{SR}.  

\begin{lem}
\label{local1}
Let $ \mathcal{X}$ and $ \mathcal{Y}$ be  subsets of $ {\mathbb Z}  /  U  {\mathbb Z} $ and  $t$ an integer satisfying $1 \leq t \leq \min_{p \in J} p^{1/2}$. Then we have that

\begin{equation}
\label{local2}
T_{J}(\mathcal{X}, \mathcal{Y}) \leq | \mathcal{X} | | \mathcal{Y} | \, \,  \exp \left( - \sum_{p \in J} \frac{1}{p} \right) \exp \left( \frac{L(\mathcal{X}, \mathcal{Y})} {t} + t \, w(J) \right) \,,
\end{equation}

\noindent
where $L(\mathcal{X}, \mathcal{Y}) = \log {( \frac{U^2}{|\mathcal{X}||\mathcal{Y}|})}$ and $w(J) = \sum_{p \in J} \frac{1}{p^2}$ .

\end{lem}

\vspace{2mm}
\noindent
Lemma \ref{local1} easily leads to the proposition below, which gives an upper bound for $T(\mathcal{X}, \mathcal{Y})$, the number of pairs $(x,y)$ in $\mathcal{X} \times \mathcal{Y}$ such that $x+y$ is an invertible element modulo $U$, that is, $x+y \not\equiv 0 \, (\text{mod} \, p)$ for all $p$ in $I$.

\begin{prop}
\label{local3}
Let $ \mathcal{X}$ and $ \mathcal{Y}$ be subsets of  $ {\mathbb Z}  /  U  {\mathbb Z} $ with $ |\mathcal{X}|, |\mathcal{Y}| \, \geq \frac{U}{M_{R}}$ . Then we have 
\begin{equation}
\label{local4}
T(\mathcal{X}, \mathcal{Y}) \leq \, \frac{\phi(P)}{P} \, |\mathcal{X}||\mathcal{Y}| \, \exp \left( \frac{36}{\log \log R} \right) \, , 
\end{equation}

\noindent
where $P=  \prod_{\substack{Q^2 < p \leq R}} \, p$ \, .
\end{prop}

\vspace{2mm}
\noindent
{\sc Proof.---}  Let $J$ be the subset of $I$ consisting of primes p such that $Q^{2} < p \leq R $. Then $$ T(\mathcal{X}, \mathcal{Y}) \, \leq \, T_{J}(\mathcal{X}, \mathcal{Y}) \, $$

\noindent
and by Lemma \ref{local1} applied to bound  $T_{J}(\mathcal{X}, \mathcal{Y})$,  we see that for any integer $1 \leq t \leq Q$ we have 

\begin{equation}
\label{local5}
T(\mathcal{X}, \mathcal{Y}) \leq | \mathcal{X} | | \mathcal{Y} | \, \,  \exp \left( - \sum_{Q^2 < p \leq R} \frac{1}{p} \right) \exp \left( \frac{L(\mathcal{X}, \mathcal{Y})} {t} + t \, w(J) \right) \,,
\end{equation}

\noindent
where $L(\mathcal{X}, \mathcal{Y}) \leq  9 \, \log R $ and $w(J) = \sum_{Q^2 < p \leq R} \frac{1}{p^2} \leq  \,  \frac{2}{Q^2}$. Thus

\begin{equation}
\label{local6}
\exp \left( - \sum_{Q^2 < p \leq R} \frac{1}{p} \right) \leq \exp \left(  \sum_{Q^2 < p \leq R} \frac{2}{p^2}  \right) \prod_{Q^2 < p \leq R} (1-p^{-1}) \leq \frac{\phi(P)}{P} \exp \left( \frac{4}{Q^2} \right) \, ,
\end{equation}

\noindent
on using $- \log(1-x) \leq x +2 x^2$, valid for $0 \leq x \leq 1/2$ . From \eqref{local5}, \eqref{local6} and the bounds for $L(\mathcal{X}, \mathcal{Y})$ and $w(J)$ we conclude that, for any integer $ 1 \leq t \leq Q $, we have

\begin{equation}
\label{local7}
T(\mathcal{X}, \mathcal{Y}) \leq \, \frac{\phi(P)}{P} \, |\mathcal{X}||\mathcal{Y}| \, \exp \left( \frac{9 \log R}{t} + \frac{9 t}{Q^2} \right) \, , 
\end{equation}

\noindent
from which \eqref{local4} follows on taking the integer $t$ in the interval $[Q/2, Q]$ and recalling that $Q = \log R \log \log R$.

\subsection{Proof of Theorem \ref{bsrthm}} \label{brsproof}

In this subsection we give a simple proof of the Theorem \ref{bsrthm}. To this end, for any $ A, B \subset{ \mathbb{N}} $ and an integer $n$ we set    

\begin{equation} \label{addirpre} 
r_{A,B}(n) = | \{(a,b)\in A \times B : a+b = n \} |,
\end{equation}
which for brevity we denote by $r(n)$. Clearly, we have that 

\begin{equation} \label{addfunctionbound}
r(n) \leq \min ( |A|, |B|),
\end{equation}

\vspace{2mm}
\noindent
and that 
\begin{equation}
\label{ot}
 P_{N;A,B} = \sum_{ \substack { 1 \leq p \leq 2N, \\ p \text { prime}}}  r(p).
\end{equation}

\vspace{2mm}
\noindent
Using the bound \eqref{addfunctionbound} for $ r(n) $ and the Chebyshev bound for the number of primes not exceeding $x$, we then get from \eqref{ot} that 

\begin{equation}
\label{trivial}
 P_{N;A,B} \ll \frac{ N } {\log{N} } \min (|A|, |B|) \ll \frac{ N }{ \log {N} } |A|^{\frac{1}{2}}|B|^{\frac{1}{2}} \ll \frac{|A||B|}{\log N} \, R \,, 
\end{equation}  

\noindent
as required.

\section{Reduction to well distributed subsets}
\label{reduct}
Let $N$, $A$ and $B$ be as in Theorem \ref{mall1}. In what follows we take  $R = \frac{1000N}{|A|^{1/2}|B|^{1/2}}$, as in Section 1, and for this $R$, we define $U, M_{R}$ and $Q$ by \eqref{umq}. Also, for any subset $Z$ of ${\mathbb Z}$, we denote by $\tilde{Z}$ the image of $Z$ in $ {\mathbb Z}/ U \mathbb{Z}$ under the natural projection map from ${\mathbb Z}$. 

\vspace{2mm}
\noindent
Let $A_{1} = A \cap [1, N^{1/8}] , A_{2}= A \cap [N^{1/8}, N]$. Then we have

\begin{equation}
\label{smbig}
P_{N;A,B} = P_{N;A_{1},B} + P_{N;A_{2},B} \leq |A_1||B| +  P_{N;A_{2},B} \ll  \,\frac{|A||B|}{\log {N}}  \, +  \, P_{N;A_{2},B},
\end{equation}

\noindent
when $N$ is large enough, since $|A_1| \leq N^{\frac{1}{8}}$, $|B| \, \leq N$ and  $ \frac{N^2}{(\log N)^2}\ll |A||B| $. 

\vspace{2mm}
\noindent
We now estimate $P_{N;A_{2},B}$. To this end, for any $a$ in  $ {\mathbb Z}/ U \mathbb{Z}$ we define 
 $m(a)$ and $n(a)$ to be, respectively,  $| \lbrace x \in A_{2} : x \, \equiv \, a  \, {\rm mod} \,  U \rbrace |$ and $| \lbrace y \in B : y \, \equiv \, a  \, {\rm mod} \,  U \rbrace |$ and then set

$$ \mathcal{C}_{1} = \lbrace a \in \tilde{A_{2}} : m(a) \leq \frac{|A_{2}|}{U} \, M_{R} \rbrace \, ,$$  $$ \mathcal{D}_{1} = \lbrace b \in \tilde{B} : n(b) \leq \frac{|B|}{U} \, M_{R} \rbrace \, ,$$
$$ \mathcal{C}_{2} = \lbrace a \in \tilde{A_{2}} : m(a) > \frac{|A_{2}|}{U} \, M_{R} \rbrace, \, \, \,   $$  $$ \mathcal{D}_{2} = \lbrace b \in \tilde{B} : n(b) > \frac{|B|}{U} \, M_{R} \rbrace \, .$$

\noindent
Since $\sum_{a \in \tilde{A}} m(a) = |A|$ and $\sum_{b \in \tilde{B}} n(b) = |B|$, it follows that  $ |\mathcal{C}_{2}|  \, \leq \, \frac{U}{M_{R}} $ and $ \, |\mathcal{D}_{2}| \, \leq \, \frac{U}{M_{R}} $ .

\vspace{2mm}
\noindent
Let us now define 
$$ A_{3}= \lbrace x \in A_{2 }:   x \, \equiv \, a  \, ({\rm mod} \,  U) \, \,  \text{for some}  \, \, a \in \mathcal{C}_{1} \rbrace \, ,$$
$$ B_{1} = \lbrace y \in B:   y \, \equiv \, b  \, ({\rm mod} \,  U) \, \,  \text{for some}  \, \, b \in \mathcal{D}_{1} \rbrace \, ,$$
$$ A_{4}= \lbrace x \in A_{2 }:   x \, \equiv \, a  \, ({\rm mod} \,  U) \, \,  \text{for some}  \, \, a \in \mathcal{C}_{2} \rbrace \, ,$$
$$ B_{2} = \lbrace y \in B:   y \, \equiv \, b  \, ({\rm mod} \,  U) \, \,  \text{for some}  \, \, b \in \mathcal{D}_{2} \rbrace \, .$$

\noindent
Then we have
\begin{equation}
\label{splitbig}
P_{N;A_{2},B} = P_{N;A_{3},B_{1}} + P_{N;A_{3},B_{2}} + P_{N;A_{4},B_{1}} + P_{N;A_{4},B_{2}} \, .
\end{equation}

\noindent
We first estimate $P_{N;A_{3},B_{2}} $. To do this, for any $a \in \tilde{A}_{3} = {\mathcal C}_1$ we define $A_{3,a}$ by $$ A_{3,a}= \lbrace x \in A_{3 }:   x \, \equiv \, a  \, ({\rm mod} \,  U) \rbrace \, ,$$ and similarly for any $b \in \tilde{B}_{2} ={\mathcal D}_2 $ we define $B_{2,b}$ by $$ B_{2,b}= \lbrace y \in B_{2 }:   y \, \equiv \, b  \, ({\rm mod} \,  U) \rbrace \, .$$

\noindent
Then we have a partition of $A_{3}$ and $B_{2}$ as follows: $$ A_{3} = \cup_{a \in \mathcal{C}_{1}}  \, A_{3,a} \,\,  \text{and} \,\, B_{2} = \cup_{b \in \mathcal{D}_{2}}  \, B_{2,b}  \,.$$

\noindent
Clearly, we then have 

\begin{equation}
\label{splitsum}
P_{N;A_{3},B_{2}} = \sum_{\substack{a \in \mathcal{C}_{1}, b \in \mathcal{D}_{2}}} P_{N;A_{3,a},B_{2,b}}  \, .
\end{equation}

\noindent
The summand on the right of \eqref{splitsum} can be estimated as

\begin{equation}
\label{splitsum1} 
P_{N;A_{3,a},B_{2,b}} \ll \frac{N}{\phi(U) \, \log N } |A_{3,a}|^{1/2}\, |B_{2,b}|^{1/2} \, ,
\end{equation}

\noindent
Indeed, if a pair $(x,y) \in A_{3,a} \times B_{2,b}$ is such that $x+y$ is a prime $p_{x,y}$, then $p_{x,y} \, \equiv \, a+b \, {\rm mod} \, U$. Since  under the condition $|A||B| \, \geq \frac{5000 N^2}{(\log N)^2}$ we have $R \leq \frac{2}{5} \log N$ and thus $U \leq N^{1/2}$ from \eqref{umq} and the Chebyshev bounds, the Brun-Titchmarsh inequality \eqref{siegel} shows that there are at most $ \frac{4N}{\phi(U) \log N}$ such primes $p_{x,y}$. Further, each such  prime can be written in at most $ \min( |A_{3,a}|, |B_{2,b}|) \leq  |A_{3,a}|^{1/2}\, |B_{2,b}|^{1/2} $ many ways as a sum $x+y$, with $x \in A_{3,a} , y \in B_{2,b}$.  
These remarks yield \eqref{splitsum1}.

\vspace{2mm}
\noindent
Using \eqref{splitsum1} in \eqref{splitsum} we then get  

\begin{equation}
\label{splitbound}
P_{N;A_{3},B_{2}}  \, \ll \, \frac{N}{\phi(U) \, \log N } \, \sum_{\substack{a \in \mathcal{C}_{1}}, b \in \mathcal{D}_{2}}  \, |A_{3,a}|^{1/2}\, |B_{2,b}|^{1/2} \, .
\end{equation}

\vspace{2mm}
\noindent
By the Cauchy-Schwarz inequality applied to the sum on the right hand side of \eqref{splitbound}, we have 

\begin{equation}
\label{cs1}
\sum_{\substack{a \in \mathcal{C}_{1}}, b \in \mathcal{D}_{2}}  \, |A_{3,a}|^{1/2}\, |B_{2,b}|^{1/2} \leq |\mathcal{C}_{1}|^{1/2} |\mathcal{D}_{2}|^{1/2}\left(\sum_{\substack{a \in \tilde{A_3}, b \in \tilde{B_2}}}  \, |A_{3,a}|\, |B_{2,b}| \right)^{1/2} 
\end{equation}

\noindent
from which and \eqref{splitbound} we get $$ P_{N;A_{3},B_{2}}  \, \ll \, \frac{N}{\phi(U) \, \log N } \, \, |\mathcal{C}_{1}|^{1/2} |\mathcal{D}_{2}|^{1/2} |A_{3}|^{1/2} |B_{2}|^{1/2} \, .$$ Now using $|\mathcal{C}_{1}| \leq U \, , |\mathcal{D}_{2}| \leq \frac{U}{M_{R}}$ and $\phi(U) \gg \frac{U}{\log R}$, which follows from Mertens formula, we get that 

$$ P_{N;A_{3},B_{2}}  \ll \frac{N}{\log N} \frac{\log R}{\sqrt{M_{R}}} \, |A|^{1/2} |B|^{1/2}  \, ,$$

\noindent
since $A_{3} \subset A , B_{2} \subset B $. Now recalling the definitions $R$ and $M_{R}$, the latter from \eqref{umq}, we get that 
$$ P_{N;A_{3},B_{2}} \ll \frac{|A||B|}{\log N}  \,  \log \log R  \, .$$

\noindent
Similarly, we get the bounds 
$$P_{N;A_{4},B_{1}} ,  P_{N;A_{4},B_{2}} \ll \frac{|A||B|}{\log N}  \,  \log \log R  \, .$$

\vspace{2mm}
\noindent
Using these bounds in \eqref{splitbig}, we then see that

\begin{equation}
\label{welldistbou}
P_{N;A_{2},B} \, \ll  \, P_{N;A_{3},B_{1}}  + \, \frac{|A||B|}{\log N}  \,  \log \log R  \, .
\end{equation}

\noindent
Thus, from \eqref{welldistbou} and \eqref{smbig},  we conclude that

\begin{equation}
\label{welldistbou1}
P_{N;A,B} \, \ll  \, P_{N;A_{3},B_{1}}  + \, \frac{|A||B|}{\log N}  \,  \log \log R  \, .
\end{equation}

\vspace{2mm}
\noindent
Therefore, to complete the  proof of Theorem \ref{mall1} we need to show that

\begin{equation}
\label{final}
P_{N;A_{3},B_{1}} \ll  \frac{|A||B|}{\log N}  \,  \log \log R  \, . 
\end{equation}

\noindent
To do this, we may assume that 

\vspace{2mm}
\noindent
(i) $|A_{3}| \geq \frac{2|A_{2}|}{U} M_{R} $ and $|B_{1}| \geq \frac{2|B|}{U} M_{R}$. 

\vspace{2mm}
\noindent
Indeed, if (i) does not hold, say, $|A_{3}| < \frac{2|A_{2}|}{U} M_{R}$, then using \eqref{trivial} and \eqref{umq} we have the stronger conclusion that  

\begin{equation}
\label{trivial1}
 P_{N;A_3,B_1} \ll  \frac{ N }{ \log {N} } |A_3|^{\frac{1}{2}}|B_1|^{\frac{1}{2}} \ll \frac{N}{\log N} \frac{|A_2|^{\frac{1}{2}}|B_1|^{\frac{1}{2}} M_{R}^{1/2}}{U^{1/2}} \ll \frac{|A||B|}{\log N} \,  \,, 
\end{equation} 

\noindent
since $U \geq e^{\frac{R}{2}}$, by the Chebyshev bound. A similar argument disposes the case $|B_{1}| < \frac{2|B|}{U} M_{R}$ as well. By the definitions of $A_3$ and $B_1$ given at the beginning of this section, we also have

\vspace{2mm}
\noindent
(ii) $ | \lbrace x \in A_{3 }:   x \, \equiv \, a  \, ({\rm mod} \,  U) \rbrace |  \, \leq \, \frac{|A_{2}|}{U} \, M_{R}$ and  $ | \lbrace y \in B_{1 }:   y \, \equiv \, b  \, ({\rm mod} \,  U) \rbrace |  \, \leq \, \frac{|B|}{U} \, M_{R}$ for any $a \in \tilde{A}_{3}, b \in \tilde{B}_{1}$ ,

\vspace{2mm}
\noindent
(iii)~ each element of $A_{3}$ is larger than $N^{1/8}$.

\vspace{2mm}
\noindent
These remarks bring us to our final section, where we shall prove \eqref{final}    taking account of the conditions (i), (ii) and (iii) above and thereby complete the proof of the Theorem \ref{mall1}.

\section{Proof Theorem \ref{mall1}}
\label{mainproof}

\vspace{2mm}
\noindent
We shall prove \eqref{final} by closely following the method of \cite{SR}. We shall assume throughout that $N$ is a sufficiently large integer. We begin by noting that if $(a,b) \in A_{3} \times B_1$ is such that $a+b$ is a prime number, then 
by (iii) above $N^{1/8} \leq a+b$ and consequently $\frac{1}{8} \log N \leq \Lambda(a+b)$, where $\Lambda$ is the Von Mangoldt function. It then follows that  

\begin{equation}
\label{denbound}
 P_{N;A_{3},B_{1}}  \, \log N   \leq \;8 \sum_{a \in A_{3}, b \in B_{1}} \Lambda (a+b) \, ,
\end{equation} 
\noindent
We have the identity $$\Lambda(n) = - \sum_{\substack{ d|n}} \mu  (d) \log{d}  \, .$$

\noindent
We now estimate right hand side  of \eqref{denbound}. To this end, we set $L=N^{1/2}$ and write $ \Lambda(n) = \Lambda^{\sharp}(n) + \Lambda^{\flat}(n)$, where  $$\Lambda^{\sharp}(n) = - \sum_{\substack{ d|n, \\ d \leq L}} \mu  (d) \log{d} \, , \quad   \Lambda^{\flat}(n) = - \sum_{\substack{ d|n,\\ d > L}} \mu  (d) \log{d} \, .$$

\noindent
Substituting $ \Lambda(n) = \Lambda^{\sharp}(n) + \Lambda^{\flat}(n)$ into \eqref{denbound}, we get that

\begin{equation}
\label{repbound}
 P_{N;A_{3},B_{1}}  \, \log N   \, \ll  \,  \sum_{n} \, r(n) \, \Lambda^{\sharp}(n) \, +  \, \sum_{n} \, r(n) \, \Lambda^{\flat}(n) \, ,
\end{equation} 
\noindent
where $r(n)$ is the number of pairs $(a,b) \in A_{3} \times B_{1}$ such that $n=a+b$. 

\vspace{2mm}
\noindent
Let us first estimate the second sum on the right of above inequality. Since $r(n)=0$ for $n$ not in the interval $[1,2N]$ , we have that 

\begin{equation}
\label{orthoind}
\sum_{n} \, r(n) \, \Lambda^{\flat}(n) \, =  \, \int_{0}^{1} \left( \sum_{n} r(n) e(-nt) \right)  \,  \left( \sum_{\substack{1 \leq n \leq 2N}} \Lambda^{\flat}(n) e(nt) \right) \,  dt  \, ,
\end{equation} 
\noindent
by orthogonality of the functions $t \mapsto e(nt)$ on $[0,1]$ . By Lemma \ref{ftbound}, we have that

\begin{equation}
\label{ftbound1}
\sum_{\substack{1 \leq n \leq 2N}} \Lambda^{\flat}(n) \, e(nt)  \, \ll  \, \frac{N}{(\log N)^{100}} \, .
\end{equation}

\noindent
From definition of $r(n)$ it immediately follows that $ \sum_{n} r(n) \, e(-nt) = \widehat{A_{3}}(-t) \widehat{B_{1}}(-t) $ . Consequently, we have from \eqref{orthoind} and \eqref{ftbound1} that

\begin{equation}
\sum_{n} \, r(n) \, \Lambda^{\flat}(n)  \, \ll \, \frac{N}{(\log N)^{100}} \, \int_{0}^{1} \, | \widehat{A_{3}}(-t) | \, |\widehat{B_{1}}(-t)| \, dt \, .
\end{equation}

\noindent
Applying Cauchy-Schwarz inequality to the above integral and using the Parseval relation together with the fact that $|A_{3}| \leq |A| $ and $ |B_{1}| \leq |B|$ we get that $$\sum_{n} \, r(n) \, \Lambda^{\flat}(n)  \, \ll \, \frac{N}{(\log N)^{100}} \, |A|^{1/2} |B|^{1/2} \,. $$ On recalling $R= \frac{1000N}{|A|^{1/2}|B|^{1/2}}$ and using the lower bound $|A| \, |B| \gg \frac{N^2}{(\log N)^2}$ , we obtain that 

\begin{equation}
\label{seconbound}
\sum_{n} \, r(n) \, \Lambda^{\flat}(n)  \, \ll \, |A| \, |B| \, .
\end{equation}

\vspace{2mm}
\noindent
Now we estimate first term on the right of \eqref{denbound}. On recalling the definition of $\Lambda^{\sharp}(n)$ we obtain

\begin{equation} 
\label{mainbound}
\sum_{n} \, r(n) \, \Lambda^{\sharp}(n) = - \sum_{ \substack{1 \leq d \leq L }} \mu(d) \, \log{d}  \sum_{ \substack{n \, \equiv \, 0  \, {\rm mod} \,d}} r(n)  \, ,
\end{equation} 

\noindent
after an interchange of summations. We note that 
 
\begin{equation} 
\label{countingmultiples}
 \sum_{ \substack{n \, \equiv \, 0 \, {\rm mod} \, d}} r(n) \, = \, \frac{1}{d} \sum_{ \substack{a \,{\rm mod} \,d }}  \, \sum_{n}  \, r(n) \, e(an/d) \, =  \, \frac{1}{d} \, \sum_{ q|d} \, \sum_{ \substack{ a \, {\rm mod^{*}} \, q} }  \, \sum_{n} r(n) \, e (an/q) \, ,  
 \end{equation}
by orthogonality of characters on the group $ \mathbb{Z} / d \mathbb{Z} $. On combining \eqref{countingmultiples} with \eqref{mainbound}, interchanging summations and recalling definition of $ \omega (q,L) $ from  \eqref{arith}, we deduce that 

\begin{equation} \label{maintermwql}
\sum_{n} r(n) \Lambda^{\sharp}(n) =   \sum_{ \substack{1 \leq q \leq L}} \omega(q,L) \sum_{ \substack{a \, {\rm mod^{*}} \, q}} \widehat{A_{3}}(a/q) \,  \widehat{B_{1}} (a/q)  \, .
\end{equation}

\noindent
We estimate the contribution to the sum on the right-hand side of \eqref{maintermwql} from $ q $ satisfying $ N^{1/8} < q \leq L $ by showing that 

\begin{equation} 
\label{bigvalueqcontribution}
 \sum_{ \substack{N^{1/8} < q \leq L}} \omega(q,L) \sum_{ \substack{a \, {\rm mod^{*}} \, q}} \widehat{A_{3}} (a/q) \widehat{B_{1}} (a/q)  \,  \ll \,  N^{7/8} \, (\log{N})^{2} \, |A|^{1/2} \, |B|^{1/2} \, \ll |A| \,|B| \, .
\end{equation}
Indeed, by \eqref{asymptotic2} we have that the absolute value of the left side of \eqref{bigvalueqcontribution} does not exceed

\begin{equation} 
\label{boundingwql}
 \frac{(\log{2L})^{2}}{N^{1/8}} \sum_{ \substack{1 \leq q \leq L}}  \, \sum_{ \substack{ a \, {\rm mod^{*}} \, q}}  \, | \widehat{A_{3}}(a/q)| \,  |\widehat{B_{1}} (a/q)|  \, \, \leq  \, \frac{(\log{2L})^{2}(N + L^2) |A_{3}|^{1/2}|B_{1}|^{1/2}}{N^{1/8}}  \, , 
\end{equation}

\noindent
where we have applied the Cauchy–Schwarz inequality followed by the large sieve  inequality \eqref{largeineq} to left-hand side of above relation. Since $ L =N^{1/2} $, we see using $|A_{3}| \leq |A|, |B_{1}| \leq |B|$  and   $|A||B| \, \gg \, \frac{N^2}{(\log N)^2}$ that \eqref{bigvalueqcontribution} follows from\eqref{boundingwql}.

\vspace{2mm}
\noindent
We now consider the contribution to the sum on the right-hand side of \eqref{maintermwql} from $ q $ in the range $ 1 \leq q \leq N^{1/8} $. We set 

\begin{equation} 
\label{ValueT}
T = \sum_{ \substack{1 \leq q \leq N^{1/8}}} \, \omega(q,L) \, \sum_{ \substack{ a \, {\rm mod^{*}} \, q }}  \widehat{A_{3}} (a/q) \, \widehat{B_{1}} (a/q) \, .
\end{equation}

\noindent
and use asymptotic formula for  $\omega(q,L)$ given by \eqref{asymptotic1}.

\noindent
The contribution of error term of this asymptotic formula for $\omega(q,L)$ to $T$ is 

\begin{equation}
\label{41} 
 \ll \frac{1}{(\log N)^{100}} \, \sum_{ \substack{1 \leq q \leq N^{1/8} }} \frac{2^{\nu(q)} \log{2q}}{q}  \, \sum_{ \substack{ a \, {\rm mod^{*}} \, q}} |  \widehat{A_{3}} (a/q)| \, | \widehat{B_{1}} (a/q)|  \, ,
\end{equation} 
 
\noindent
where we have used the value $L=N^{1/2}$ and $\alpha = 100$. Using the trivial bound $ 2^{ \nu(q)} \log{2q} \, \ll  q $  we see that \eqref{41} is 
$$ \ll \frac{1}{(\log N)^{100}} \, \sum_{ \substack{1 \leq q \leq N^{1/8} }} \, \sum_{ \substack{ a \, {\rm mod^{*}} \, q}} |  \widehat{A_{3}} (a/q)| \, | \widehat{B_{1}} (a/q)| \ll \frac{N}{(\log N)^{100}} |A|^{1/2} \, |B|^{1/2} \ll |A| |B|  ,$$

\noindent
where we have used the Cauchy-Schwarz inequality, large sieve inequality \eqref{largeineq} and the bound $|A||B| \gg \frac{N^2}{(\log N)^2}$ . Thus, we have 

\begin{equation}
\label{asympT}
T = \sum_{ \substack{1 \leq q \leq N^{1/8}}} \, \frac{\mu(q)}{\phi(q)} \, \sum_{ \substack{ a \, {\rm mod^{*}} \, q }}  \widehat{A_{3}} (a/q) \, \widehat{B_{1}} (a/q) + O \left( |A||B| \right) \, .
\end{equation}

\vspace{2mm}
\noindent
We recall from \eqref{umq} that $U=\prod_{p \leq R} \, p$ . Then since $R=\frac{1000 \,N}{|A|^{1/2}|B|^{1/2}} $ and $|A||B| \gg \frac{N^2}{(\log N)^2}$, we see that $U \leq N^{1/8}$ for all large $N$. We set $T(U)$ to be the sum over $ q $ on the right-hand side of \eqref{asympT} restricted to $q | U $. Since for all other $ q $ we have either $ \mu(q) = 0 $ or $ q > R $, the triangle inequality applied to \eqref{asympT} shows that

\begin{equation}
\label{differTU}
|T-T(U)| \, \ll \, \sum_{ \substack{R \leq q \leq N^{1/8}}} \, \frac{1}{\phi(q)} \, \sum_{ \substack{ a \, {\rm mod^{*}} \, q }}  |\widehat{A_{3}} (a/q)| \, |\widehat{B_{1}} (a/q)|  \, +  \, O \left( |A||B| \right) \, .
\end{equation}

\noindent
We shall estimate the sum over $ q $ in \eqref{differTU} by using $ \frac{q}{\log{\log{q}}} \ll \phi(q) \ll q $, the Cauchy-Schwarz inequality and the large sieve inequality \eqref{largeineq}. Since $ \frac{\log{\log{q}}}{q} $ is decreases with $ q $ for $  q \geq 10 $, we get that 
\begin{align}
|T-T(U)| & \ll \frac{\log \log R}{R} \, N \, |A|^{1/2} \, |B|^{1/2} \\
         & \ll |A| \, |B| \, \log \log R  \label{difftut} \, ,   
\end{align}

\noindent
since $R=\frac{1000 \,N}{|A|^{1/2}|B|^{1/2}}$.

\vspace{2mm}
\noindent
Now we estimate $T(U)$. A simple argument using standard properties of Ramanujan sums, given below (3.23) on page 969 of \cite{SR} shows  that

\begin{equation} 
\label{coprimeeletoU}
T(U) = \frac{U}{\phi(U)} |\{ (a,b) \in A_{3} \times B_{1} : (a+b, U) = 1 \}| \, .
\end{equation}

\vspace{2mm}
\noindent
As before, we use $\tilde{A}_{3}$ to denote the image of $A_{3}$ under the natural projection from the set of all integers ${\mathbb Z}$ to ${\mathbb Z} / U {\mathbb Z}$  and similarly denote by $\tilde{B}_{1}$ the image of $B_{1}$ . Further, for any residue class $a$ modulo $U$, let $m_{A_{3}}(a)$ be the number of elements of the set $A_{3}$ that belongs to this residue class. Similarly, we define $m_{B_{1}}(b)$ for any residue class $b$ modulo $U$ . Let $D_{1} = \frac{|A_{2}|}{U} \, M_{R}$ and $D_{2} = \frac{|B|}{U} \, M_{R}$ . We then have using condition (ii) given at the end of the preceding section that 

$$ \sum_{a \in \tilde{A}_{3}} \, m_{A_{3}} (a) = |A_{3}|  \, \, \text{with} \, \,  0 \leq m_{A_{3}}(a) \leq D_{1}  \, ,$$ and
$$ \sum_{b \in \tilde{B}_{1}} \, m_{B_{1}} (b) = |B_{1}|  \, \, \text{with} \, \,  0 \leq m_{B_{1}}(b) \leq D_{2}  \, .$$

\vspace{2mm}
\noindent
Let us set $c(a,b)$ to be $1$ when $a+b$ is invertible modulo in ${\mathbb Z} / U {\mathbb Z}$ and to be $0$ otherwise. Then from \eqref{coprimeeletoU} we get that 

\begin{equation}
\label{localT}
T(U) =  \, \frac{U}{\phi(U)} \, \sum_{(a,b) \in \tilde{A}_{3} \times \tilde{B}_{1}}  c(a,b) \,  m_{A_{3}}(a) \, m_{B_{1}}(b) \, .
\end{equation}

\noindent
We estimate above sum with the help of the optimization principle given in Subsection \ref{opt}. From the Lemma \ref{bil}, we then have that 

\begin{equation}
\label{optimabound}
\sum_{(a,b) \in \tilde{A}_{3} \times \tilde{B}_{1}}  c(a,b) \,  m_{A_{3}}(a) \, m_{B_{1}}(b) \, \leq \, \sum_{(a,b) \in \tilde{A}_{3} \times \tilde{B}_{1}}  c(a,b) \, x_{a}^{*} \, y_{b}^{*}  \, ,
\end{equation}

\noindent
for some $x_{a}^{*}$ and $y_{b}^{*}$ with $a$ varying over $\tilde{A}_{3}$ and $b$ varying over $\tilde{B}_{1}$ , satisfying the following conditions. All the $x_{a}^{*}$ are either $0$ are $D_{1}$, excepting at most one, which must lie in $(0,D_{1})$  and similarly, all $y_{b}^{*}$ are either $0$ or $D_{2}$, excepting at most one, which must lie in $(0,D_{2})$ . Moreover, if $\mathcal{X}$ and $\mathcal{Y}$ denote,  respectively, the subsets of $\tilde{A}_{3}$ and $\tilde{B}_{1}$ for which $x_{a}^{*} \neq 0$ and $y_{b}^{*} \neq 0$ , then $ |\mathcal{X}| \, D_{1} \geq |A_{3}| \, \geq \, (|\mathcal{X}|-1) \,  D_{1}$ and $ |\mathcal{Y}| \, D_{2}  \, \geq |B_{1}| \, \geq  \,(|\mathcal{Y}|-1) \, D_{2}$ . Thus, from this we have the bounds 
$$ D_{1} \, \leq \, \frac{A_{3}}{|\mathcal{X}|-1} \, \leq \, \frac{2 \, |A_{3}|}{|\mathcal{X}|}$$ and $$ D_{2} \, \leq \, \frac{B_{1}}{|\mathcal{Y}|-1} \, \leq \, \frac{2 \, |B_{1}|}{|\mathcal{Y}|} \, ,$$

\noindent
where we use $|\mathcal{X}| \geq \frac{|A_{3}|}{D_{1}} \geq 2 $ and $|\mathcal{Y}| \geq \frac{|B_{1}|}{D_{2}} \geq 2 $, valid by condition (i) given at the end of Section \ref{reduct}. These bounds on $D_{1}, D_{2}$ together with \eqref{localT} and \eqref{optimabound} give

\begin{equation}
T(U) \, \ll \, \frac{U}{\phi(U)} \, \frac{|A_{3}||B_{1}|}{|\mathcal{X}||\mathcal{Y}|} \, \sum_{(a,b) \in \mathcal{X} \times \mathcal{Y}} c(a,b) \, ,
\end{equation}

\vspace{2mm}
\noindent
Note that $\mathcal{X},\mathcal{Y}$ are subsets of ${\mathbb Z} / U {\mathbb Z}$ with $|\mathcal{X}| \geq \frac{|A_{3}|}{D_{1}} \geq \frac{U}{M_{R}} $ and $|\mathcal{Y}| \geq \frac{|B_{1}|}{D_{2}} \geq \frac{U}{M_{R}} $ and that the sum on the right of above relation is nothing but $T(\mathcal{X},\mathcal{Y})$ of Proposition \ref{local3}, which gives  

\begin{equation}
\label{tubound}
T(U) \, \ll \, \frac{U}{\phi(U)} \, \frac{\phi(P)}{P} \, |A_{3}| \, |B_{1}|  \, \exp \left( \frac{36}{\log \log R} \right) \, ,
\end{equation}

\noindent
where $P =\prod_{Q^2 < p \leq R} \, , U=\prod_{p \leq R}$ and $Q=\log R \, \log \log R$. By Merten's formula we then get 

\begin{equation}
\label{finaltubound}
T(U) \, \ll \, |A| \, |B|  \,  \log \log R \, \exp \left( \frac{C}{\log \log R} \right)  \ll \, |A| \, |B| \, \log \log R  \,.
\end{equation}

\noindent
From \eqref{finaltubound}, \eqref{difftut}, \eqref{bigvalueqcontribution}, \eqref{seconbound} and \eqref{repbound} we then conclude that
$$ P_{N;A_{3},B_{1}} \, \ll \, \frac{|A| \,|B|}{\log N} \, \log \log R \, ,$$ which together with \eqref{welldistbou1} yields Theorem \ref{mall1}.

\vspace{3mm}
\noindent
{\bf Acknowledgement :} I am thankful to Prof. D. S. Ramana and Prof. Gyan Prakash for their guidance and suggestions during this work and also for their comments on the final draft of this paper. I am very grateful to the referee for carefully reading this article and providing me with an extensive list of suggestions. I also wish to thank the Harish-Chandra Research Institute for the excellent facilities given to me.

\vspace{5mm}
\Address

\end{document}